\documentclass[a4paper,11pt]{amsart}
\usepackage{amssymb,amsfonts,amsmath}
\usepackage{layout}
\usepackage[latin1]{inputenc}

\def\m{\mathcal}

\def\C{\mathbb{C}}
\def\c2{\mathbb{C}^2}
\def\R{\mathbb{R}}
\def\N{\mathbb{N}}

\def\N{\mathbb{N}}
\def\P{\mathbb{P}}

\def\1{\bold{1}}

\def\a{\alpha}

\def\e{\varepsilon}

\def\f{\varphi}

\def\p{\psi}

\def\om{\omega}

\newcommand \W {\Omega}

\newcommand \mE {\mathcal E}

\newcommand \vphi {\varphi}

\def\E{{\mathcal{E}_{\chi}(X,\om)}}

\newtheorem{lem}{Lemma}[section]
\newtheorem{pro}[lem]{Proposition}

\newtheorem{def/not}[lem]{Definition/Notations}
\numberwithin{equation}{section}
\newtheorem{thm}[lem]{Theorem}
\newtheorem{cor}[lem]{Corollary}

\newtheorem{exa}[lem]{Example}

\newenvironment{proof3.1}
{\noindent {\it{Proof of theorem 3.1}}}{$\Box$ \linebreak[4]}

\begin{document}

\title[A priori estimates for solutions of Monge-Amp\`ere equations]
{A priori estimates for weak solutions of complex Monge-Amp\`ere equations}

\author{S.BENELKOURCHI \& V.GUEDJ \& A.ZERIAHI}

\begin{abstract}
Let $X$ be a compact K\"ahler manifold and $\om$ a smooth closed form of bidegree
$(1,1)$ which is nonnegative and big. We study
the classes ${\mathcal E}_{\chi}(X,\om)$
of $\om$-plurisubharmonic functions of finite weighted Monge-Amp\`ere energy.
When the weight $\chi$ has fast growth at infinity, the corresponding functions 
are close to be bounded.

We show that if a positive Radon measure is suitably dominated by the Monge-Amp\`ere
capacity, then it belongs to the range of the Monge-Amp\`ere operator on some
class ${\mathcal E}_{\chi}(X,\om)$. 
This is done by establishing a priori estimates on the
capacity of sublevel sets of the solutions.

Our result extends U.Cegrell's and S.Kolodziej's results
and puts them into a unifying frame. It also gives a simple proof of
S.T.Yau's celebrated a priori ${\mathcal C}^0$-estimate.
\end{abstract}

\maketitle


\section{Introduction} 
Let $X$ be a compact connected  K\"ahler manifold of dimension $n \in \N^* $.
Throughout the article $\om$ denotes a smooth closed form of bidegree
$(1,1)$ which is nonnegative and {\it big}, i.e. such that
$\int_X \om^n >0$. We continue the study started in
 \cite{GZ 2}, \cite{EGZ} of the complex Monge-Amp\`ere equation
 $$
\! \! \!\!  \! \! \!\!  \! \! \!\!  \! \! \!\!  \! \! \!\!  
\! \! \!\!  \! \! \!\!  \! \! \!\!  \! \! \!\!  \! \! \!\!  
\! \! \!\!  \! \! \!\!  \! \! \!\!  \! \! \!\!  \! \! \!\!  
\! \! \!\!  \! \! \!\!  \! \! \!\!  
\mbox{(MA)}_\mu \hskip3cm (\om +dd^c \f )^n  = \mu ,
$$
where $\f$, the unknown function, is $\om$-plurisubharmonic: 
this means that $\f \in L^1(X)$ is upper semi-continuous and $\om + dd^c \f \ge 0$ 
is a positive current. We let $PSH(X, \om )$ denote the set of all such functions
 (see \cite{GZ 1} for their basic properties). 
Here $\mu $ is a fixed positive Radon measure of total mass
$\mu(X) = \int _X\om ^n $, and $d=\partial+\overline{\partial}$,
$d^c=\frac{1}{2i\pi}(\partial-\overline{\partial})$. 

Following \cite{GZ 2} we say that a $\om$-plurisubharmonic function
$\f$  has finite weighted Monge-Ampère energy, 
$\f \in {\mathcal E}(X,\om)$, when its Monge-Ampère measure
$(\om + dd^c\f )^n $ is well defined, and there exists 
an increasing function  $\chi :\ \R^- \to \R^- $ 
such that $\chi ( -\infty) = -\infty $  and $ \chi \circ \f \in L^1((\om + dd^c \f )^n) $. 
In general $\chi $ has very slow growth at infinity, so that $\f $ is far from 
being bounded.

The purpose of this article is twofold. First we extend one of the main results
of \cite{GZ 2} by showing
\vskip.2cm

\noindent {\bf THEOREM A.} 
{\it There exists $\f \in {\mathcal E}(X,\om)$ such that
$\mu=(\om+dd^c \f)^n$ if and only if $\mu$ does not charge pluripolar sets.
}
\vskip.2cm

\noindent This results has been established in [GZ 2] when $\om$ is a K\"ahler form. It is important for applications to complex dynamics and K\"ahler geometry to consider as well forms $\om$ that are less positive (see [EGZ]).

We then look for conditions on the measure $\mu$ which
 insure that the solution $\f$ is almost bounded. 
Following the seminal work of S. Kolodziej [K 2,3], we say that $\mu $ is dominated by 
the Monge-Ampère Capacity $Cap _{\om}$ if there exists a function 
$ F : \R^+ \to \R^+ $ such that $\lim_{t \to 0^+} F(t)=0$ and
\begin{equation} \tag{\dag}
\mu(K ) \le F(Cap_{\om} (K) ), \quad \mbox{for all Borel subsets } \ K \subset X.
\end{equation}
Here $ Cap _\om $ denotes the global version of the
Monge-Ampère capacity introduced by E.Bedford and
A.Taylor \cite{BT} (see section 2). 

Observe that $\mu $ does not charge
pluripolar sets since $F(0)=0.$  When $F(x) \lesssim x^\alpha $ vanishes 
at order 
$\alpha >1$ and $ \om$ is K\"ahler, S. Kolodziej has proved \cite{K 2} that the solution 
$\f \in PSH(X, \om ) $ of (MA)$_\mu$ is {\it continuous}. 
The boundedness part of this result was extended in \cite{EGZ} to the case when $\om$ is merely big and 
nonnegative. 
If $F(x) \lesssim x^\alpha $  with $0< \alpha <1 ,$ two of us have proved in
 \cite{GZ 2} that the solution $\f $ has finite $\chi -$energy, where 
$\chi (t) = -(-t)^p , \ p=p(\alpha )>0$.
This result was first established by U. Cegrell in a local context \cite{Ce}.

Another objective of this article is to fill in the gap inbetween Cegrell's and 
Kolodziej's results, by considering all intermediate dominating functions $F.$ 
Write $ F_\e (x) = x [\e( -\ln (x) /n)]^n$ where  $\e:\R \rightarrow  [0 ,\infty [$ 
is nonincreasing. Our second main result is:
\vskip.2cm

\noindent {\bf THEOREM B.} 
{\it If  $\mu(K) \leq F_{\e}(Cap_{\om}(K))$
for all Borel  subsets $K  \subset X$,
then $\mu=(\om+dd^c \f)^n$ where $\f \in PSH(X,\om)$ satisfies  
 $\sup _X \f = 0$ and
\begin{equation*}
Cap_{\om}(\f<-s) \leq  \exp (-nH^{-1}(s)).
\end{equation*}
Here $H^{-1}$ is the reciprocal function of $H(x) = 
e \int_{0}^x \e (t) dt + s_0 ,$  
where $s_0 = s_0 ( \e , \om) \ge 0 $ only depends on $\e $ and $\om .$
}
\vskip.2cm

This general statement has several useful consequences:
\begin{itemize}
\item if  $ \int_0^{+\infty } \e (t) dt <+\infty , $ then $H^{-1}(s)= +\infty $
for $s \ge s_\infty := e \int_{0}^{+\infty } \e (t) dt+s_0,$ hence
$Cap _\om (\f < -s ) =0.$ 
This means that $\f $ is bounded from below by $-s_{\infty} .$ 
This result is due to S. Kolodziej [K 2,3] when $\om$ is K\"ahler,
and \cite{EGZ} when $\om \geq 0$ is merely big;

\item the condition (\dag) is easy to check for measures with density
in $L^p$, $p>1$. 
Our  result thus gives a simple proof (Corollary 3.2), following the seminal approach of S. Kolodziej ([K2]),
of the $\m C^0 $-a priori estimate of 
S.T. Yau \cite{Y}, which is crucial for proving the Calabi conjecture 
(see \cite{T} for an overview);

\item when  $ \int_0^{+\infty } \e (t) dt = +\infty ,$
 the solution $\f $ is generally unbounded. 
The faster $\e(t)$ decreases towards zero,
 the faster the growth of $H^{-1}$ at infinity, 
hence the closer is $\f $ from being bounded;

\item the special case $\e \equiv 1$ is of particular interest. 
Here  $\mu (\cdot) \le Cap_\om  (\cdot)$, 
and our result shows that $Cap_\om (\f <-s )$ decreases exponentially fast, 
hence $\f $ has `` loglog-singularities''. These are the type of singularities of the metrics used in 
Arakelov geometry in relation with measures $\mu = f d V $  whose density has Poincaré-type singularities 
(see [Ku], [BKK]).
\end{itemize}

We prove Theorem B in {\it section 3}, after establishing Theorem A
in {\it section 2.1} and recalling some useful
facts from [GZ 2], [EGZ] in {\it section 2.2}.
We then test the sharpness of our estimates in {\it section 4}, where we give 
examples of measures 
fulfilling our assumptions: these are absolutely continuous with respect to $\om^n$,
and their density do not belong to $L^p$, for any $p>1$.

\section{Weakly singular quasiplurisubharmonic functions}

The class ${\mathcal E}(X,\om)$ of $\om$-psh functions with finite
weighted Monge-Amp\`ere energy has been introduced and studied in
\cite{GZ 2}. It is the largest subclass of
$PSH(X,\om)$ on which the complex Monge-Amp\`ere operator
$(\om+dd^c \cdot)^n$ is well-defined and the comparison principle
is valid. Recall that $\f \in {\mathcal E}(X,\om)$
if and only if $(\om+dd^c \f_j)^n(\f \leq -j) \rightarrow 0$,
where $\f_j:=\max(\f,-j)$.

\subsection{The range of the Monge-Amp\`ere operator}

The range of the operator $(\om+dd^c \cdot)^n$ acting on
${\mathcal E}(X,\om)$ has been characterized in [GZ 2] when $\om$
is a {\it K\"ahler} form. We extend here this result to the case when
$\om$ is merely nonnegative and big.

\begin{thm}
Assume $\om$ is a smooth closed nonnegative (1,1) form on $X$,
and $\mu$ is a positive Radon measure such that
$\mu(X)=\int_X \om^n>0$.

Then there exists $\f \in {\mathcal E}(X,\om)$ such that
$\mu=(\om+dd^c \f)^n$ if and only if $\mu$ does not charge 
pluripolar sets.
\end{thm}

\begin{proof}
We can assume without loss of generality that $\mu $ and $\om $ are normalized so that 
$\mu (X) = \int_X\om ^n =1.$ Consider, for $A>0$,
$$
\m C _A(\om ): = \{ \nu\  \text{probability measure } / \ \nu (K) \le A \cdot Cap_\om (K), 
\text{for all }\ K\subset X \},
$$
where $Cap_\om $ denotes the Monge-Ampère capacity introduced by E.Bedford and A.Taylor
in [BT] (see \cite{GZ 1} for this compact setting). Recall that
$$
Cap _\om (K) : = \sup \left\{ \int _K (\om + dd^c u )^n \ / \  u\in PSH(X, \om), \ 0
\le u \le 1 \right\}.
$$

We first show that a measure $\nu \in \m C _A(\om ) $ is the Monge-Ampère of a  function
$ \psi \in \m E^p (X , \om ) , $ for any $0<p <1$, where 
$$
\m E^p (X , \om) := \{ \psi \in \mE(X, \om)\ / \  \psi \in L^p \big ( (\om + dd^c \p)^n \big )\}.
$$

Indeed, fix $\nu \in \m C _A(\om ), \ 0<p<1,$ and $\om_j : = \om + \e_j \W $, 
where $\W$  is a k\"ahler form on $X$, and $\e_j>0 $ decreases towards zero. 
Observe that 
$PSH(X, \om ) \subset  PSH(X, \om_j ),$ hence $ Cap _\om (.) \le  Cap _{\om_j} (.),$ so that $\nu \in \m C _A(\om_j ).$ It follows from Proposition 3.6 and 2.7 in \cite{GZ 1} that there exists $C_0 >0$ such that for any $v \in  PSH(X, \om_j )$ normalized by $\sup _X v = -1,$ we have 
$$
Cap _{\om_j} (v<-t) \le \frac{C_0}{t}, \  \text{for all } t\ge 1.
$$
This yields $ \m E^p (X , \om_j) \subset L^p (\nu )$: if $v\in   \m E^p (X , \om_j)$
 with  $\sup _X v = -1,$ then
\begin{eqnarray*}
\int _X (-v)^p d \nu  &=& p\cdot \int_0 ^{+\infty} t^{p-1} \nu (v < -t)dt \\ 
&\le & p A \cdot \int_1 ^{+\infty} t^{p-1} Cap _\om  (v < -t)dt  + C_p\\ 
& \le & \frac{pA C_0}{1-p} + C_p <+\infty.
\end{eqnarray*}
It follows therefore from Theorem 4.2 in \cite{GZ 2} that there exists 
$\vphi _j \in  \m E^p (X , \om_j) $  with $\sup_X \vphi _j =-1$ and 
$(\om_j +dd^c \vphi_j)^n = c_j \cdot \nu ,$ where $c_j = \int_X\om_j ^n \ge 1$ 
decreases towards 1 as $\e_j $  decreases towards zero. 
We can assume without loss of generality that $1\le c_j \le 2.$ Observe that the $\vphi _j$'s 
have uniformly bounded energies, namely
$$
\int_X (-\f_j)^p (\om_j +dd^c \vphi_j)^n \le 2 \int_X (-\f_j)^p d\nu 
\le 2 \left[\frac{pAC_0}{1-p} +C_p \right].
$$
Since $\sup_X\f_j =-1 ,$ we can assume (after extracting a convergent subsequence) that 
$\f_j \to \f $ in $ L^1(X),$ where $\f \in PSH(X, \om ), \ \sup_X \f =-1.$ 

Set $\phi _j : = (\sup_{l\ge j} \f_l)^* .$ 
Thus $\phi_j \in PSH(X, \om_j),$ and $\phi_j $ decreases towards $\f .$
Since $\phi_j \ge \f_j,$ it follows from the ``fundamental inequality'' (Lemma 2.3 in \cite{GZ 2}) that 
$$
\int_X (-\phi_j)^p (\om_j +dd^c \phi_j)^n \le 2^n \int_X (-\f_j)^p (\om_j +dd^c \f_j)^n
\le C^\prime <+\infty.
$$
Hence it follows from stability properties of the class $\m E^p (X , \om) $ that
$\f \in  \m E^p (X , \om)$ (see Proposition 5.6 in \cite{GZ 2}). Moreover 
$$
(\om_j +dd^c \phi_j)^n \ge \inf_{l\ge j}(\om_l +dd^c \f_l)^n \ge \nu ,
$$
hence $(\om +dd^c \vphi)^n = \lim (\om_j +dd^c \phi_j)^n \ge \nu.$ 
Since $\int_X \om^n = \nu(X) = 1,$ this yields $\nu = (\om +dd^c \vphi)^n$ as claimed above.

We can now prove the statement of the theorem. 
One implication is obvious: if 
$\mu =(\om +dd^c \vphi)^n, \ \f \in \mE( X, \om),$ then $\mu $ does not charge pluripolar sets, 
as follows from Theorem 1.3 in \cite{GZ 2}. 

So we assume now $\mu$ that does not charge pluripolar sets. 
Since $\m C _1(\om ) $ is a compact convex set of probability measures which contains 
all measures $(\om +dd^c u)^n$, $ u\in PSH(X, \om ), \ 0\le u\le 1,$ we can project 
$\mu $ onto $\m C_1(\om ) $ and get, by a generalization of Radon-Nikodym theorem (see [R], [Ce]), 
$$
\mu = f\cdot \nu , \ \nu \in \m C_1(\om ) , \ 0\le f \in L^1(\nu ).
$$
Now $ \nu = (\om +dd^c \psi)^n $ for some $\psi \in  \m E^{1/2} (X , \om), \ \psi \le 0 ,$ 
as follows from the discussion above. Replacing $\psi$ by $e^\psi $ shows that
 we can actually assume  $\psi $ to be bounded (see Lemma 4.5 in [GZ 2]). 
We can now apply line by line the same proof as that of Theorem 4.6 in \cite{GZ 2} 
to conclude that $\mu = (\om + dd^c \f )^n $  for some $\f \in \mE (X, \om ).$
\end{proof}

\subsection{High energy and capacity estimates}

Given $\chi:\R^- \rightarrow \R^-$ an increasing function, we consider,
following \cite{GZ 2},
$$
\E:=\left\{ \f \in {\mathcal E}(X,\om) \, / \,
\int_X (-\chi) (-|\f|) \, (\om+dd^c \f)^n <+\infty \right\},
$$
 Alternatively a function
$\f \leq 0$ belongs to $\E$ if and only if
$$
\sup_j \int_X (-\chi) \circ \f_j \, (\om+dd^c \f_j)^n <+\infty,
\text{ where } \f_j:=\max(\f,-j)
$$
is the {\it canonical approximation} of $\f$
by bounded $\om$-psh functions.
When $\chi(t)=-(-t)^p$, $\E$ is the class ${\mathcal E}^p(X,\om)$
used in previous section.

The properties of classes $\E$ are quite different whether
the weight $\chi$ is convex (slow growth at infinity)
or concave. In previous works [GZ 2], 
two of us were mainly interested in weights $\chi$
of moderate growth at infinity (at most polynomial).
Our main objective in the sequel is to construct solutions
$\f$ of $(MA)_{\mu}$ which are ``almost bounded'', i.e. in
classes $\E$ for concave weights $\chi$ of arbitrarily high growth.

For this purpose it is useful to relate the property $\f \in \E$ to the
speed of decreasing of $Cap_{\om}(\f<-t)$, as $t \rightarrow +\infty$.
We set
$$
\hat{{\mathcal E}}_{\chi}(X,\om):=\left\{ \f \in PSH(X,\om) \, / \,
\int_0^{+\infty} t^n \chi'(-t) Cap_{\om}(\f<-t) dt<+\infty \right\}.
$$

An important tool in the study of classes $\E$ are the ``fundamental inequalities''
(Lemmas 2.3 and 3.5 in [GZ 2]), which allow to compare the weighted energy
of two $\om$-psh functions $\f \leq \p$. These inequalities are only valid for
weights of slow growth (at most polynomial), while they become immediate
for classes $\hat{{\mathcal E}}_{\chi}(X,\om)$. So are the convexity properties
of $\hat{{\mathcal E}}_{\chi}(X,\om)$.
We summarize this and compare these classes in the following:

\begin{pro}
The classes $\hat{{\mathcal E}}_{\chi}(X,\om)$ are convex and stable under
maximum: if $\hat{{\mathcal E}}_{\chi}(X,\om) \ni \f \leq \p \in PSH(X,\om)$,
then $\p \in \hat{{\mathcal E}}_{\chi}(X,\om)$.

One always has $\hat{{\mathcal E}}_{\chi}(X,\om) \subset \E$, while
$$
{\mathcal E}_{\hat{\chi}}(X,\om) \subset \hat{{\mathcal E}}_{\chi}(X,\om),
\text{ where } \chi'(t-1)=t^n \hat{\chi}'(t).
$$
\end{pro}

Since we are mainly interested in the sequel in weights with 
(super) fast growth at infinity,
the previous proposition shows that $\hat{{\mathcal E}}_{\chi}(X,\om)$ and $\E$
are roughly the same: a function $\f \in PSH(X,\om)$
belongs to one of these classes if and only if $Cap_{\om}(\f<-t)$ decreases
fast enough, as $t \rightarrow +\infty$.

\begin{proof}
The convexity of $\hat{{\mathcal E}}_{\chi}(X,\om)$ follows from
the following simple observation: if $\f,\p \in \hat{{\mathcal E}}_{\chi}(X,\om)$
and $0 \leq a \leq 1$, then
$$
\left\{ a\f+(1-a)\p <-t \right\} \subset
\left\{ \f<-t \right\} \cup  \left\{ \p <-t \right\}.
$$
The stability under maximum is obvious.

Assume $\f \in \hat{{\mathcal E}}_{\chi}(X,\om)$.
We can assume without loss of generality $\f \leq 0$ and $\chi(0)=0$.
Set $\f_j:=\max(\f,-j)$. It follows from Lemma 2.3 below that
\begin{eqnarray*}
\int_X (-\chi) \circ \f_j \, (\om+dd^c \f_j)^n &=&
\int_0^{+\infty} \chi'(-t) (\om+dd^c\f_j)^n(\f_j < -t) dt \\
&\leq& \int_0^{+\infty} \chi'(-t) t^n Cap_{\om}(\f<-t) dt
<+\infty, 
\end{eqnarray*}

This shows that $\f \in \E$.
The other inclusion goes similarly, using the second inequality
in Lemma 2.3 below.
\end{proof}

If $\f \in \E$ (or $\hat{{\mathcal E}}_{\chi}(X,\om)$), then the bigger the growth
of $\chi$ at $-\infty$, the smaller $Cap_{\om}(\f<-t)$ when $t\rightarrow +\infty$, hence
the closer $\f$ is from being bounded. Indeed $\f \in PSH(X,\om)$
is bounded iff it belongs to $\E$ for all weights $\chi$,
as was observed in [GZ 2], Proposition 3.1. Similarly
$$
PSH(X,\om) \cap L^{\infty}(X)=\bigcap_{\chi} \hat{{\mathcal E}}_{\chi}(X,\om),
$$
where the intersection runs over all concave increasing functions
$\chi$.
\vskip.2cm

We will make constant use of the following result:

\begin{lem}
Fix $\f \in {\mathcal E}(X,\om)$. Then  for all $s>0$ and $0 \leq t \leq 1$,
$$
t^n Cap_{\om}(\f<-s-t) \leq \int_{(\f<-s)} (\om+dd^c \f)^n \leq s^n Cap_{\om}(\f<-s),
$$
where the second inequality is true only for $s \geq 1$.
\label{ine}
\end{lem}

The proof is a direct consequence of the comparison principle (see Lemma 2.2 in [EGZ] and [GZ 2]).

\section{Measures dominated by capacity}

 From now on $\mu$ denotes a positive Radon measure on $X$ whose total mass is $Vol_\om (X)$:
 this is an obvious necessary condition
in order to solve $(MA)_{\mu}$. To simplify numerical computations, we assume in the sequel
that $\mu $ and $\om $ have been normalized so that $$\mu(X) =Vol_\om (X)  =\int_X \om^n =1. $$

When $\mu=e^h \om^n$ is a smooth volume form
and $\om$ is a K\"ahler form, S.T.Yau has proved [Y]
that $(MA)_{\mu}$ admits a unique {\it smooth} solution
$\f \in PSH(X,\om)$ with $\sup_X \f=0$.
Smooth measures are easily seen to be nicely dominated by
the Monge-Amp\`ere capacity (see the proof of Corollary 3.2 below).

\vskip.1cm

Measures dominated by the Monge-Amp\`ere capacity have been extensively studied by
S.Kolodziej in [K 2,3,4]. Following S. Kolodziej ([K3], [K4]) with slightly different notations, fix $\e:\R \rightarrow  [0 ,
\infty [$ a continuous decreasing  function and set 

$$F_{\e}(x):=x [\e(-\ln x/n)]^{n}, x > 0.$$
We will consider probability measures $\mu$ satisfying the following condition :
for all Borel subsets $K \subset X$,
$$
\mu(K) \leq F_{\e}(Cap_{\om}(K)).$$
 The main result  achieved
in [K 2], can be formulated as follows:  
If $\om$ is a K\"ahler form and $\int_0^{+\infty} {\e(t)}dt  <+\infty$
then $\mu=(\om+dd^c \f)^n$ for some {\it continuous} function
$\f \in PSH(X,\om)$.

\vskip.1cm

The condition $\int_0^{+\infty} {\e(t)}dt  <+\infty$ means
that $\e$ decreases fast enough towards zero at infinity. This gives
a quantitative estimate on how fast $\e( -\ln Cap_{\om}(K)/n)$,
hence $\mu(K)$, decreases towards zero as $Cap_{\om}(K) \rightarrow 0$.

When $\int_0^{+\infty} \e(t)dt=+\infty$, it follows from Theorem 2.1
that $\mu=(\om+dd^c \f)^n$ for some function $\f \in {\mathcal E}(X,\om)$, 
but $\f$ will generally be unbounded. Our second main result
measures how far $\f$ is from being bounded:

\begin{thm}\label{subest}
Assume for all compact subsets $K  \subset X$,
\begin{equation}
\mu(K) \leq F_{\e}(Cap_{\om}(K)).
\label{dom}
\end{equation}
Then $\mu=(\om+dd^c \f)^n$ where $\f \in {\mathcal E}(X,\om)$ is 
such that $\sup _X \f = 0$ and
$$
Cap_{\om}(\f<-s) \leq  \exp (-nH^{-1}(s)),
\text{ for all } s>0.
$$
Here $H^{-1}$ is the reciprocal function of
 $H(x) =e \int_{0}^x  \e(t) dt + s_0$,
 where 
$s_0 = s_0 ( \e , \om ) \ge 0 $
 is a constant which only depends on $\e $ and $\om .$ 

In particular $\f \in \E$ where
$-\chi(-t)=\exp (  n H^{-1}(t)/2)$.
\end{thm}

\noindent Recall that here, and troughout the article, $\om \geq 0$ is merely big.

Before proving this result we make a few observations.
\begin{itemize}
\item It is interesting to consider as well the case when  $\e(t)$  increases towards
 $+\infty$.  One can 
 then obtain solutions $\f$ such that
$Cap_{\om}(\f<-t)$ decreases at a polynomial rate. 
When e.g. $\om$ is K\"ahler and $\mu(K) \leq Cap_\om(K)^{\a}$,
$0<\a<1$, it follows from Proposition 5.3
in [GZ 2] that $\mu=(\om+dd^c \f)^n$ where $\f \in {\mathcal E}^p(X,\om)$
for some $p=p_{\a}>0$. Here  $ {\mathcal E}^p(X,\om)$ denotes the Cegrell type class 
 $ {\mathcal E}_\chi (X,\om),$ with $\chi(t) = -(-t)^p.$

\item When $\e(t) \equiv 1$, $F_{\e}(x)= x$ and $H(x) \asymp  e.x$.
Thus Theorem \ref{subest} reads 
$ \mu \leq Cap_{\om} \Rightarrow \mu=(\om+dd^c \f)^n$, where
$$
Cap_{\om}(\f<-s) \lesssim \exp\left( -n s/e \right).
$$
This is precisely the rate of decreasing corresponding to functions
which look locally like $-\log(-\log||z||)$, in some local
chart $z \in U \subset \C^n$.
This class of $\om$-psh functions with ``loglog-singularities''
is important for applications (see [Ku], [BKK]).

\item If $\e(t)$ decreases towards zero, then
$Cap_{\om}(\f<-t)$ decreases at a superexponential rate.
The faster $\e(t)$ decreases towards zero, the slower the growth
of $H$, hence the faster the growth of $H^{-1}$ at infinity.
When $\int^{+\infty} \e(t)dt<+\infty$, the function $\e$ decreases so fast
that $Cap_{\om}(\f<-t)=0$ for $t>>1$, thus $\f$ is bounded.
This is the case when $\mu(K) \leq Cap_{\om}(K)^{\a}$ for some
$\a>1$ \cite{K 2}, \cite{EGZ}.

\item When $\int^{+\infty} \e(t)dt=+\infty$, the solution $\f$ may well be unbounded
(see Examples in section 4).
At the critical case where $\mu \leq F_{\e}(Cap_{\om})$
for all functions $\e$ such that $\int^{+\infty} \e(t)  dt=+\infty$,
we obtain
$$
\mu=(\om+dd^c \f)^n
\text{ with }
\f \in PSH(X,\om) \cap L^{\infty}(X),
$$
as follows from Proposition 3.1 in [GZ 2].
This partially explains the difficulty in describing the range
of Monge-Amp\`ere operators on the set of {\it bounded}
(quasi-)psh functions.
\end{itemize}

\begin{proof} The assumption on $\mu $ implies in particular that it
 vanishes on pluripolar sets. It follows from Theorem 2.1  that there exists
 a function $\f \in \m E ( X, \om )$ such that $\mu=(\om+dd^c \f)^n$
and $\sup _X \f =0.$
Set
 $$
g (s) : = -\frac{1}{n} \log Cap_{\om}(\f<-s),  \ \ \forall s>0
.$$
The function $g$ is increasing on $[0 , +\infty ] $ and
$ g (+\infty)=+\infty$, since $Cap_{\om}$ vanishes on pluripolar sets.
Observe also that $g (s) \geq 0$ for all $s \geq 0$, since
$$
g (0)=-\frac{1}{n} \log Cap_{\om}(X)=-\frac{1}{n} \log Vol_{\om}(X)=0.
$$

It follows from  Lemma \ref{ine} and (\ref{dom}) that for all $s>0$ and $0 \leq t \leq 1$,
$$
t^n Cap_{\om}(\f<-s-t) \leq \mu(\f<-s) \leq F_{\e}\left(Cap_{\om}(\f<-s)\right).
$$
Therefore for all $s>0$ and $0 \leq t \leq 1$,
\begin{equation}
\log t - \log \e \circ g (s) + g (s) \leq g (s+t).
\label{est}
\end{equation}

We define an increasing sequence $(s_j)_{j \in \N}$ by induction
setting 
$$ 
s_{j+1} = s_{j} + e \e \circ g (s_{j}), \text{ for all } j \in \N.
$$
\vskip.2cm

\noindent {\it The choice of $s_0$}. Recall that (\ref{est}) is only valid for
 $0\le t\le 1.$ 
We choose $s_0\ge 0 $ large enough so that
 \begin{equation}
 e .  \e \circ g (s_{0}) \le 1.
\end{equation} 
This will allow us 
to use (\ref{est}) with $t = t_j =  s_{j+1} -s_j \in [0 , 1]$, since  
$ \e \circ g $ is decreasing, while $s_j \ge s_0 $ is increasing, hence 
$$
0 \le t_j =  e \e \circ g (s_{j})  \le  e \e \circ g (s_{0}) \le  1.
$$
We must insure that $s_0 = s_0 ( \e , \om ) $ can chosen to be independent of $\f.$ 
This is a consequence of Proposition 2.7 in \cite{GZ 1}: since $\sup _X \f = 0,$ 
there exists $c_1 (\om ) >0 $ so that 
$0 \le \int_X (-\f) \om ^n \le c_1 (\om ),$ hence 
$$
g (s) : = -\frac{1}{n} \log Cap_{\om}(\f<-s) \ge \frac{1}{n} \log s  - \frac{1}{n}\log ( n + c_1 (\om )).
$$
Therefore $g (s_0 ) \ge \e ^{-1}(1/e)$ for $s_0 = s_0(\e , \om ) :=  (n + c_1 (\om)) \exp (n \e^{- 1} (1 \slash e)),
$ 
which is independent of $\f$. This  yields 
 $ e .  \e \circ g (s_{0}) \le 1$, as desired.

\vskip.2cm 
\noindent {\it The growth of $s_j$. } We can now apply (\ref{est}) and get
 $g (s_j) \ge j + g (s_0) \ge j.$
Thus  $\lim g (s_j)=+\infty$.
There are two cases to be considered. 

If $s_{\infty}=\lim s_j \in \R^+$, then
$g (s) \equiv +\infty$ for $s > s_\infty$, i.e. 
$ Cap_{\om}(\f<-s)=0,  \ \ \forall s  > s_\infty$. 
Therefore $\f$ is bounded from below by $-s_\infty$,
in particular $\f \in \mE_\chi(X,\om)$ for all $\chi .$

Assume now (second case) that $s_j\to +\infty.$ For each   
 $s>0, $  there exists $ N = N_s  \in \N $ such that
 $s_{N} \le s < s_{N +1}.$
We can estimate $s \mapsto N_s$:
\begin{eqnarray*}
  s \le s_{N +1}&=&  \sum_0^{N} (s_{j+1} -s_j) + s_0  = \sum_{j = 0}^{N}
   e \,  \e \circ g (s_{j}) + s_0 \\
&  \le &  e \sum_{0}^{N}
 \e (j) + s_0 
\le   e. \e(0)  + e \int_{0}^{N}
  \e (t)dt + s_0  =:   H(N) ,
\end{eqnarray*}
Therefore $H^{-1}(s) \le N \le g (s_N) \le g (s),$ hence
\begin{equation*}
Cap_{\om}(\f<-s) \leq  \exp (-nH^{-1}(s)).
\end{equation*}

Set now $-\chi(-t) =\exp (  n  H^{-1}(t)/2)$. Then
\begin{multline*}
\int_0^{+\infty} t^n \chi'(-t) Cap_{\om}(\f<-t) dt\\  \le
\frac{n}{2}\int_0^{+\infty} {  t^n}\frac{1}{\e(  H^{-1}(t)) + \tilde{s}_0 } \exp (-n H^{-1}(t)/2) dt\\
\le C \int_{0}^{+\infty}  { t^n} \exp (-nt/2 ) dt <+\infty.
\end{multline*}
This shows that  $\f \in \E$ where $\chi(t)=-\exp (  n  H^{-1}(-t)/2)$.

It follows from the proof above that when $ \int_0^{+ \infty} \e (t) d t < + \infty$, the solution $\f$ is bounded since in this case  we have
$$
s_{\infty} := \lim_{j \to + \infty} s_j \leq   s_0 (\e,\om) + e \ \e (0) + e  \int_0^{+ \infty} \e (t) d t   < + \infty
$$
 where $s_0 (\e,\om)$ is an absolute constant satisfying $(3.3)$ (see above). 
\end{proof}

Let us emphasize that Theorem 3.1 also yields a slightly simplified proof of 
the following result [K 2], [EGZ]:
if $\mu(K) \leq F_{\e}(Cap_{\om}(K))$ for some decreasing function $\e:\R \rightarrow \R^+$ 
such that $\int^{+\infty} \e(t) dt <+\infty$, then the sequence $(s_j)$ above
is convergent, hence $\mu=(\om+dd^c \f)^n$, where $\f \in PSH(X,\om)$
is {\it bounded}. 
For the reader's  convenience we indicate a  proof of the following 
important particular case:

\begin{cor}
Let $\mu=f \om^n$ be a measure with density $0 \leq f \in L^p(\om ^n )$, 
where $p>1$ and $\int_X f \om^n=\int_X \om^n$.
Then there exists a unique bounded function $\f \in PSH(X,\om)$ such that
$(\om+dd^c \f)^n=\mu$, $\sup_X \f=0$ and
$$
0\le ||\f||_{L^\infty (X)}  \leq  C(p,\om) . ||f||^{1/n}_{L^p(\om ^n)} ,
$$
where $C(p,\om)>0$ only depends on  $p $ and $\om$.
\end{cor}

This a priori bound is a crucial step in the proof by S.T.Yau of the Calabi conjecture
(see [Ca], [Y], [A], [T], [Bl]). The proof presented here follows Kolodziej's new and decisive pluripotential approach (see [K2]).
Let us stress that the dependence $\om \longmapsto C (p,\om)$ is quite explicit, as we shall see in the proof.
This is important when considering degenerate situations  [EGZ].

\begin{proof} 
We claim that there exists $C_1(\om)$ such that
\begin{equation}
\mu(K) \leq \Big [C_1(\om)||f||^{1/n}_{L^p(\om ^n )}\Big ]^n \left[ Cap_{\om} (K) \right]^2,
\; \text{ for all Borel sets } K \subset X. 
\end{equation}
Assuming this for the moment, we can apply Theorem 3.1 with
 $\e(x)=C_1(\om)||f||^{1/n}_{L^p(\om ^n )}
 \exp(-x)$,
which yields, as observed at the end of the proof of Theorem 3.1  
$$
||\f||_{L^\infty (X)}  \leq M (f,\om) ,
$$
where $M (f,\om) := s_0 (\e,\om) + e \ \e (0) + e  \int_0^{+ \infty} \e (t) d t = s_0 (\e,\om) + 2 e C_1(\om)||f||^{1/n}_{L^p(\om ^n )}$
and  $s_0 = s_0 (\e,\om)$
is a large number $s_0 > 1$ satisfying the inequality $(3.3)$.

In order to give the precise dependence of the uniform bound $M (f,\om)$ on the $L^p-$norm of the density $f$, we need to choose $s_0$ more carefully. 
Observe that condition $(3.3)$ can be written
$$Cap_{\om} (\{\f \leq - s_0\}) \leq \exp (- n \e^{- 1} (1 \slash e).$$ 
Since $ n \e^{- 1} (1 \slash e)  =  \log \Bigl(e^n C_1 (\om)^n \Vert f \Vert_{L^p(\om^n)}\Bigr),$ we must choose $s_0 > 0$ so that
\begin{equation}
Cap_{\om} (\{\f \leq - s_0\}) \leq \frac{1}{e^n C_1 (\om)^n \Vert f \Vert_{L^p(\om^n)}}. 
\end{equation}

 We claim that for any $N \geq 1$ there exists a uniform constant $ C_2 (N,p,\om) > 0$ such that for any $s > 0,$

\begin{equation}
 Cap_{\om} (\{\f \leq - s\}) \leq C_2 (N,p) \ s^{- N} \ \Vert f\Vert_{L^p(\om^n)}. 
\end{equation}
Indeed observe first that by H\"older inequality, 
$$
\int_X (- \f)^N \om_{\f}^n =   \int_X (-\f)^N f \om^n \leq \Vert f \Vert_{L^p(\om^n)}
 \Vert \f \Vert_{L^{N q}(\om^n)}^N.
$$
Since $\f$ belongs to the compact family $\{ \psi \in PSH (X,\om) ; \sup_X \psi = 0 \}$ ([GZ2]), there exists a uniform constant $C'_2 (N,p,\om) > 0$ such that $\Vert \f \Vert_{L^{N q}(\om^n)}^N \leq C'_2 (N,p,\om)$, hence
$$ \int_X (- \f)^N \om_{\f}^n \leq C'_2 (N,p,\om) \Vert f\Vert_{L^p(\om^n)}.$$
Fix  $u \in PSH (X,\om)$ with $ - 1 \leq u \leq 0$ and  $N \geq 1$ to be specified later. If follows from Tchebysheff and energy  inequalities ([GZ2]) that
\begin{eqnarray*}
\int_{\{\f \leq - s\}} (\om + dd^c u)^n \ & \leq & \ s^{- N}  \int_X (- \f)^N (\om + dd^c u)^n \\
 & \leq &  \ c_N \ s^{- N} \max \left\{\int_X (- \f)^N \om_{\f}^n , \int_X (- u)^N \om_{u}^n \right\} \\
& \leq & \  c_N \ s^{- N} \  \max \left\{C'_2 (N,p,\om) ,1\right\}  \Vert f \Vert_{L^p(\om^n)}.
\end{eqnarray*}
We have used here the fact that $  \Vert f \Vert_{L^p (\om^n)} \geq 1$, which follows from the normalization :  $ 1 = \int_X \om^n = \int_X f \om^n  \leq \Vert f \Vert_{L^p (\om^n)}$. This proves the claim. 

Set $N = 2 n$, it follows from $(3.6)$ that  $s_0  :=  C_1 (\om)^n e^n C_2 (2 n, p,\om) \Vert f \Vert_{L^p(\om^n)}^{1 \slash n}$ satisfies the required condition $(3.5)$, which implies the estimate of the theorem.
\vskip.2 cm
We now establish the estimate $(3.4)$. Observe first that H\"older's inequality yields
\begin{equation}
\mu(K) \leq ||f||_{L^p(\om ^n )} \left[ Vol_{\om}(K) \right]^{1/q},
\text{ where } 1/p+1/q=1.
\end{equation}

Thus it suffices to estimate the volume $Vol_{\om}(K)$. Recall the definition of the Alexander-Taylor 
capacity, $T_{\om}(K):=\exp(-\sup_X V_{K,\om})$, where
$$
V_{K,\om}(x):=\sup \{ \p(x) \, / \, \p \in PSH(X,\om), \p \leq 0 \text{ on } K \}.
$$
This capacity is comparable to the Monge-Amp\`ere capacity,
as was observed
by H.Alexander and A.Taylor [AT] (see Proposition 7.1 in [GZ 1] for this compact setting):
\begin{equation}
T_{\om}(K) \leq e \exp \left[ - \frac{1}{Cap_{\om}(K)^{1/n}} \right].
\end{equation}

It thus remains to show that $Vol_{\om}(K)$ is suitably bounded
from above by $T_{\om}(K)$. This follows from Skoda's uniform integrability result:
set
$$
\nu(\om):=\sup \left\{ \nu(\p,x) \, / \, \p \in PSH(X,\om), \, x \in X \right\},
$$
where $\nu(\p,x)$ denotes the Lelong number of $\p$ at point $x$.
This actually only depends on the cohomology class $\{ \om \} \in H^{1,1}(X,\R)$.
It is a standard fact that goes back to H.Skoda (see [Z]) that there exists $C_2(\om)>0$
so that
$$
\int_X \exp\left(-\frac{1}{\nu(\om)} \p \right) \, \om^n \leq C_2(\om),
$$
for all functions $\p \in PSH(X,\om)$ normalized by $\sup_X \p=0$. We infer
\begin{equation}
Vol_{\om}(K) \leq \int_K \exp\left(-\frac{1}{\nu(\om)} V_{K,\om}^* \right) \, \om^n
\leq C_2(\om) [T_{\om}(K)]^{1/\nu(\om)}.
\end{equation}
It now follows from (3.7), (3.8), (3.9), that
$$
\mu(K) \leq ||f||_{L^p} [C_2(\om)]^{1/q} e^{1/q\nu(\om)}
\exp\left[ - \frac{1}{q\nu(\om) Cap_{\om}(K)^{1/n}} \right].
$$
The conclusion follows by observing that 
$\exp(-1/x^{1/n}) \leq C_{n} x^2$ for some explicit constant $C_{n}>0$.
\end{proof}

\section{Examples}

\subsection{Measures invariant by rotations}

In this section we produce examples of radially invariant functions/measures
which show that our previous results are essentially sharp.
The first example is due to S.Kolodziej [K 1].

\begin{exa}
We work here on the Riemann sphere $X=\P^1(\C)$, with $\om=\om_{FS}$,
the Fubini-Study volume form.
Consider $\mu=f \om$ a measure with density $f$ which is smooth and positive
on $X \setminus \{p\}$, and such that
$$
f(z) \simeq \frac{c}{|z|^2 (\log|z|)^2}, \; c>0,
$$
in a local chart near $p=0$. A simple computation yields
$\mu=\om+dd^c \f$, where $\f \in PSH(\P^1,\om)$ is smooth in $\P^1 \setminus \{p\}$
and $\f(z) \simeq -c' \log(-\log|z|)$ near $p=0$, $c'>0$,
hence
$$
\log Cap_{\om}(\f<-t) \simeq -t,
$$
Here $a \simeq b$ means that $a/b$ is bounded away from zero and infinity.

This is to be compared to our estimate 
$\log Cap_{\om}(\f<-t) \lesssim -t/e $
(Theorem \ref{subest} ) which can be applied,
as it was shown by S.Kolodziej in [K 1] that $\mu \lesssim Cap_{\om}$.
Thus Theorem \ref{subest} is essentially sharp when $\e \equiv 1$.
\end{exa}

We now generalize this example and show that the estimate provided by
Theorem \ref{subest} is
essentially sharp in all cases.

\begin{exa}
Fix $\e $ as in Theorem \ref{subest}. 
Consider $\mu =f \om$  on $X=\P^1(\C)$, where $\om=\om_{FS}$
is the Fubini-Study volume form, $f \geq 0$ is continuous on
$\P^1 \setminus \{p\}$, and 
$$
f(z)  \simeq \frac{\e  (\log (-\log |z|))}{|z|^2 (\log |z|)^2  } 
$$
in local coordinates near $p=0$.
Here $\e:\R \rightarrow \R^+$ decreases towards $0$ at $+\infty$.
We claim that there exists $A>0$ such that
\begin{equation} \label{cap}
\mu(K) \le A  {Cap_{\om}(K)}{\e (- \log Cap_{\om}(K))},
\text{ for all } K \subset X.
\end{equation}

This is clear outside a small neighborhood of $p=0$
since the measure $\mu $ is there
dominated by a smooth volume form.
So it suffices to establish this estimate when $K$ is included in a local chart near $p=0$.
Consider
$$
\tilde{K} := \{r\in [0, R]\ ;\ K \cap \{|z|=r\}\not =  \emptyset\}.
$$

It is a classical fact (see e.g. [Ra])
that the logarithmic capacity $c(K)$ of $K$  
can be estimated from below by the length of $\tilde{K}$, namely
$$
\frac{l(\tilde{K})}{4} \leq c(\tilde{K}) \leq c(K).
$$
Using that $\e$ is decreasing, hence $0 \leq -\e'$, we infer
\begin{eqnarray*}
\mu (K) &\le & 2\pi \int_0^{l(\tilde{K})} f(r) r dr
\\
&\le & {2\pi} \int _0^{ l({\tilde{K}})}\frac{\e  (\log (-\log r)) -
 \e^{\prime}(\log -\log r)}{r
 (\log r)^2  } dr \\
&=& 2\pi \frac{\e  (\log (-\log l({\tilde{K}}) )) }{-
 \log l({\tilde{K}})   }  
\le   2 \pi \frac{ \e  (\log (-\log  4 c (K)) ) }{
 -\log 4 c (K)  } .
\end{eqnarray*}

Recall now that the logarithmic capacity c(K) is equivalent to Alexander-Taylor's capacity 
$T_{\Delta} (K)$, which in turn is equivalent to the global Alexander-Taylor
capacity $T_{\om}(K)$ (see [GZ 1]):
$
c(K) \simeq T_{\Delta}(K) \simeq T_{\om}(K).
$
The Alexander-Taylor's comparison theorem \cite{AT} reads
$$
-\log 4c(K) \simeq -\log T_{\om}(K) \simeq 1/Cap_{\om}(K), 
$$
thus $\mu(K) \leq A Cap_{\om}(K) \e(-\log Cap_{\om}(K))$.
\vskip.15cm

We can therefore apply Theorem \ref{subest}. It guarantees that
$\mu=(\om+dd^c \f)$, where $\f \in PSH(\P^1,\om)$ satisfies
$\log Cap_{\om}(\f<-s) \simeq   -nH^{-1}(s)$,
with $H (s)= eA \int_{0}^s {\e(t)} dt+s_0$.
On the other hand a  simple computation shows that
$\f $ is continuous 
in $\P^1 \setminus \{p\}$ and
$$
\f  \simeq  - H(\log (- \log |z|)) \ , \ \  near \
  p=0.
$$

The sublevel set 
$(\f<-t) $ therefore coincides with the ball of  radius 
$\exp (-\exp( H^{-1}(t)))$, hence 
$\log Cap_{\om}(\f<-s) \simeq  -H^{-1}(s).$
\end{exa}

\subsection{Measures with density}
Here we consider the case when $\mu=f dV$ is absolutely continuous
with respect to a volume form.

\begin{pro}
Assume $\mu=f \om^n$ is a probability measure whose density satisfies
$f [\log(1+f)]^n \in L^1(\om^n)$. Then $\mu \lesssim Cap_{\om}$.

More generally if $f [\log(1+f)/\e(\log (1 + |\log f|))]^n \in L^1(\om^n)$
for some continuous decreasing function $\e:\R \rightarrow \R^+_*$, then
for all $K \subset X$,
$$
\mu(K) \leq F_{\e}(Cap_{\om}(K)),
\text{ where } F_{\e}(x)=Ax \left[\e\left(-\frac{\ln x}{n}\right)\right]^n, \, A>0.
$$
\end{pro}

\begin{proof}
With slightly different notations,
the proof is identical to that of Lemma 4.2 in [K 4] to which we refer the reader.
\end{proof}

We now give examples showing that Proposition 4.3 is almost optimal.

\begin{exa}
For simplicity we give local examples. The computations to follow can also
be performed in a global compact setting.

Consider $\f(z)=-\log (-\log ||z||)$, where
$||z||=\sqrt{|z_1|^2+\ldots +|z_n|^2}$ denotes the Euclidean
norm in $\C^n$. One can check that $\f$ is plurisubharmonic in
a neighborhood of the origin in $\C^n$, and that
there exists $c_n>0$ so that
$$
\mu:=(dd^c \f)^n=f \, dV_{\text{eucl}},
\text{ where }
f(z)=\frac{c_n}{||z||^{2n} (-\log||z||)^{n+1}}.
$$
 Observe that   $f [\log(1+f)]^{n-\alpha} \in L^1, \ \forall \alpha >0$
but
$f [\log(1+f)]^n \not \in L^1. $  

When $n=1$ it was observed by S. Kolodziej \cite{K 1}
that $ \mu(K) \lesssim Cap_{\om}(K).$ Proposition 4.3 yields here 
$$ 
\mu(K) \lesssim Cap_{\om}(K) (|\log Cap_{\om}(K)| + 1 ).
$$

For $n\ge 1,$ it follows from Proposition 4.3 and Theorem \ref{subest} that 
$$
\log Cap_{\om}(\f<-s) \lesssim   -nH^{-1}(s).
$$
On the other hand, one can directly check 
that $\log Cap_{\om}(\f<-s) \simeq   -nH^{-1}(s).$

One can get further examples by considering $\f(z)=\chi \circ \log||z||$, so that
$$
(dd^c \f )^n =\frac{ c_n (\chi ^{\prime}\circ \log ||z|| )^{n-1}
 \chi ^{\prime \prime}(\log ||z|| )}{||z||^{2n}}  dV_{\text{eucl}}.
$$
\end{exa}

\vskip .2cm

Slimane BENELKOURCHI \& Vincent GUEDJ \& Ahmed ZERIAHI

Laboratoire Emile Picard

UMR 5580, Universit\'e Paul Sabatier

118 route de Narbonne

31062 TOULOUSE Cedex 09 (FRANCE)

benel@math.ups-tlse.fr

guedj@math.ups-tlse.fr

zeriahi@math.ups-tlse.fr

\end{document}